\documentclass{amsart}

\usepackage{latexsym}
\usepackage{amssymb}
\usepackage{amsmath}
\usepackage{amscd}
\usepackage{amsthm}

\newtheorem{thm}{Theorem}[section]
\newtheorem{definition}[thm]{Definition}
\newtheorem{propo}[thm]{Proposition}
\newtheorem{lemma}[thm]{Lemma}
\newtheorem{corollary}[thm]{Corollary}

\newtheorem*{rac}{Remark}

\newcommand{\xteich}{\mathsf{Teich}}

\newcommand{\xml}{\mathsf{ML}}

\newcommand{\xrp}{\mathbb{RP}}
\newcommand{\xcp}{\mathbb{CP}}

\newcommand{\Mm}{\mathcal M}
\newcommand{\Ss}{\mathcal S}

\newcommand{\Nn}{\mathcal N}

\newcommand{\mr}{\mathbb{R}}

\newcommand{\mh}{\mathbb{H}}

\renewcommand{\Re}{{\mathbb R}}         
\newcommand{\la}{\langle}               
\newcommand{\ra}{\rangle}               

\newcommand{\Ric}{\text{\rm Ric}}	
\newcommand{\Ham}{\mathcal H}		
\newcommand{\tr}{\text{\rm tr}}		
\newcommand{\dev}{\mathcal D}		
		
\newcommand{\Vol}{\text{\rm Vol}}	
\newcommand{\Area}{\text{\rm Area}}     

\newcommand{\Mink}{\mathbb M}		

\newcommand{\HH}{\mathbb H}

\newcommand{\tM}{\tilde M}

\newcommand{\fin}{\hspace*{\fill} 
\quad\hbox{\hskip 1pt\vrule width 4pt height 6pt
          depth 1.5pt\hskip 1pt} \medskip }

\begin{document}

\title[Notes on a paper of Mess]
{Notes on a paper of Mess}

\author[Andersson, Barbot, et al]{Lars Andersson}
\thanks{The first author was supported in part by the NSF,
contract no. DMS 0407732.}
\address{Department of Mathematics\\
University of Miami\\
Coral Gables, FL 33124\\
USA \and 
Albert Einstein Institute \\ 
am M\"uhlenberg 1\\
D-144476 Potsdam\\
Germany}
\email{larsa@math.miami.edu}
\author[]{Thierry Barbot}
\thanks
{The second author was supported 
by CNRS, ACI ``Structures g\'eom\'etriques et Trous Noirs''.}
\address{CNRS, UMR 5669\\
Ecole Normale Sup\'erieure de Lyon\\ 
46 all\'ee d'Italie, 69364 Lyon}
\email{Thierry.Barbot@umpa.ens-lyon.fr}
\author[]{Riccardo Benedetti}
\address{Dipartimento di Matematica\\
Universit\`a di Pisa\\
Largo B. Pontecorvo, 5, I-56127 Pisa}
\email{benedett@dm.unipi.it}
\author[]{Francesco Bonsante}
\address{Dipartimento di Matematica\\
Universit\`a di Pisa\\
Largo B. Pontecorvo, 5, I-56127 Pisa}
\email{f.bonsante@sns.it}
\author[]{William M. Goldman}
\address{Department of Mathematics\\
University of Maryland\\
College Park, MD 20742}
\email{wmg@math.umd.edu}
\author[]{Fran\c{c}ois Labourie}
\address{Topologie et Dynamique\\
Universit\'e Paris-Sud\\
F-91405 Orsay (Cedex)}
\email{francois.labourie@math.u-psud.fr}
\author[]{Kevin P. Scannell}
\address{Department of Mathematics and Computer Science\\
Saint Louis University\\
St. Louis, MO 63103}
\email{scannell@slu.edu}
\author[]{Jean-Marc Schlenker}
\address{Institut de Math\'ematiques\\
Universit\'e Paul Sabatier\\
31062 Toulouse Cedex 9}
\email{schlenker@picard.ups-tlse.fr}


\date{\today}
\maketitle

\renewcommand{\thesection}{N.\arabic{section}}
\renewcommand{\thethm}{\thesection.\arabic{thm}}
\section{Introduction}\label{sec1}

In his 1990 paper ``Lorentz Spacetimes of Constant Curvature'' \cite{Me90},
Geoff Mess offered what was, at the time, a completely new approach
to the study of spacetimes in
$2+1$-dimensions, primarily by employing
tools and techniques from low-dimensional geometry and topology.
Among the many interesting results in the paper is
a complete geometric parameterization of the moduli space of
flat $2+1$ spacetimes with closed Cauchy surfaces.  
By ``geometric'' we mean that the parameterization is
defined in terms of certain objects that arise in the
Cauchy horizons of such spacetimes (specifically
$\mathbb{R}$-trees, or, dually, measured laminations -- structures
introduced by Thurston in his groundbreaking work on 
hyperbolic $3$-manifolds in the 1970's).
This allows one to do much more than merely count
degrees of freedom in the moduli space.
For example, certain previously inaccessible problems
can be resolved quite easily in this language, 
such as results on the inextendability of domains of
dependence past the Cauchy horizon (Cosmic Censorship);
see \S\ref{SS:CC} below.

The paper was written around the same time that Witten and
others began analyzing $2+1$-gravity as a kind of
toy model for quantization, and it attracted a lot of
attention (and many citations) from the physics community.
It has been similarly influential within mathematics,
with applications to $3$-dimensional affine geometry,
the deformation theory of hyperbolic manifolds, etc.

Unfortunately, the paper was never published.  In addition, 
it has never been particularly easy to find in 
preprint form since there was no electronic version available. 
For this reason, and because of the profound impact the paper has had,
we decided that it would be a useful service
to the community
to have the paper in print, even in its somewhat unpolished state.

The version of the paper that follows is nearly
identical to the original; we merely corrected
typographical errors and occasional notational 
mistakes, and also updated the references in 
the bibliography.  The numbering of definitions,
propositions, and theorems is unchanged, 
so that earlier citations to the preprint
can easily be traced in this version.

These ``Notes'' are intended as a roadmap to Mess'
paper, which can be difficult to follow at times.
We offer some simplifications to the
original arguments, an occasional 
correction, and most importantly a guide
to the literature on this subject that has appeared 
in the intervening 16 years.

\S\ref{fuchsian} - \S\ref{S:classif} of the present paper correspond
to the sections of Mess' paper with the same numbers,
for easy cross-referencing while reading the latter.
These notes also contain two additional sections 
describing some recent generalizations and
improvements upon Mess' results.
In \S\ref{S:wick} we discuss an
important unifying notion (the {\em canonical Wick rotation})
due to Benedetti and Bonsante.
The final section (\S\ref{S:CMC}) considers
slicings of spacetimes
by constant mean curvature surfaces,
a problem raised by Mess but not addressed 
directly in his paper.

Definitions, propositions, and theorems in these Notes
are numbered with a prefixed ``N''
(Lemma N.3.1, Corollary N.3.2, etc.),
to avoid confusion with the numbering in Mess' paper. 

The authors would like to thank 
Dick Canary, David Garfinkle, Steve Harris, Misha Kapovich,
Steve Kerckhoff, Cyril Lecuire, John Ratcliffe, and Ser Tan  
for helpful correspondence.
Thanks to Taejung Kim for helping typeset
Mess' manuscript, and to the referee for a number of useful
comments.

\section{Fuchsian holonomy}
\label{fuchsian}

The main result of this section is Proposition 1; this says that the
linear part of the holonomy of a closed spacelike surface with $g > 1$
in a flat $2+1$ spacetime is Fuchsian. 
The proof is based on the Milnor-Wood inequality \cite{Wo71}
and Goldman's theorem \cite{Go88} and there is a terse exposition
of these results (the proof of Goldman's theorem being essentially Matsumoto's
\cite{Ma87a}).
These are best read from the original sources.

A second proof of discreteness of the linear holonomy is sketched
at the end of \S8; see \cite{Ba05} for details.

Proposition 2 is the simple observation that any Fuchsian
group gives a flat spacetime by quotienting the interior of 
the future (or past) light cone.

\section{Realization of holonomy homomorphisms}
\label{realization}

The main result is Proposition 3, which is an application of the holonomy theorem to say that a small affine deformation of a given Fuchsian linear representation of a closed surface group is realized by a future complete flat spacetime;
rescaling then
realizes any affine deformation in this way.

The fact that the surface is closed is critical here, since there are 
examples of affine deformations of Fuchsian groups with
parabolics that do not leave invariant a future complete
regular convex domain \cite{Ba05}, as well as the ``crooked plane'' examples
of Goldman and Drumm \cite{DG90} which,
in contrast with the examples constructed by Mess,
actually act discontinuously on all of Minkowski space.

Also critical to the later development is Lemma 1, which shows that
a spacelike immersion of a surface into Minkowski space such that
the induced metric is complete must actually be an embedding.
In later sections this is applied to the developing map of the
universal cover of a Cauchy surface.

We provide here a particularly nice geometric argument generalizing
Lemma 1; this argument was known previously (at least in 
\cite{Ha87}).

This argument applies in any Lorentzian manifold $M$ admitting 
a complete Killing vector field $X$
which is everywhere timelike (in particular, non-singular), and such that, if 
$\Phi^t$ denotes the flow generated by $X$,
every orbit of $\Phi^t$ is wandering, i.e., at any point $x$ of 
$M$ there is a small
transversal $D$ to $\Phi^t$ containing $x$ and such that $D$ 
intersects every $\Phi^t$-orbit at most once. 
Observe that if 
$M$ is strongly causal (see \S\ref{s.defgh}),
then this hypothesis is automatically satisfied.

This hypothesis implies that the orbit space $Q^{\Phi}$ is a manifold,
though possibly non-Hausdorff: non-Hausdorffness might appear from \emph{ancestral pairs\/}; see \cite{HL01}.
It can be shown that under reasonable hypotheses, valid
in Minkowski and anti-de~Sitter space, $Q^\Phi$ is Hausdorff: for example, 
if $M$ admits a time function and is timelike or null geodesically complete
(Theorem $2$ in \cite{GH97}).

Let $\pi: M \rightarrow Q^\Phi$ be the projection map. For any tangent 
vector $\bar{v}$ to $Q^{\Phi}$ at a point $\theta$, define $g_0(\bar{v})$ to be the norm in $M$ of any vector $v$ orthogonal to $X(x)$
such that $x \in \theta$ and $d\pi(v) = \bar{v}$. Since $X$ is Killing, this norm does not depend on the choice of the point $x$ in the orbit $\theta$. It defines a Riemannian metric $g_0$ on $Q^\Phi$. Observe also that $\pi$ is a locally trivial fibration.
The proof of the following lemma is straightforward:

\begin{lemma}
\label{increase}
Let $S$ be a Riemannian manifold, and let 
$f: S \rightarrow M$ be a codimension $1$ isometric immersion. Then,
 the composition $\pi \circ f$ is distance increasing.
\fin
\end{lemma}

\begin{corollary}
If the Riemannian metric on $S$ is complete, then $\pi \circ f$ is a covering map by a standard argument \cite{KN63}.\fin
\end{corollary}

\begin{corollary}
\label{icon}
If $Q^{\Phi}$ is simply connected and $S$ is complete, then
$f$ is an embedding, and the image of $f$ is the graph of a section 
of $\pi$. \fin
\end{corollary}

A particularly suitable case is the \emph{static\/} case, i.e., the case where $M$ contains a hypersurface $S_0$ everywhere orthogonal to $X$ (in other words, the case where the distribution $X^\perp$ is integrable).  Then, the restriction to $S_0$ of $\pi$ is an isometric identification between $S_0$ and $Q^\Phi$. Pushing along $\Phi^t$, we obtain an identification between $M$ and the product $S_0 \times {\mathbb R}$. The vector field $X$ is then $\partial/\partial{t}$, and the Lorentzian metric at a point $(x,t)$ is:
$$g_0(x) - a(x)dt^2$$
where $a(x)$ is the opposite of the norm of $X(x,t)$. Let $g_1$ be the metric $\frac{1}{a}g_0$ on $S_0 \approx Q^{\Phi}$. Corollary~\ref{icon} can be formulated as follows:

\begin{corollary}
When $X$ is static, every complete spacelike immersed surface is the graph of a $1$-contracting map $s: (S_0, g_1) \rightarrow \mathbb R$.\fin
\end{corollary}

In Minkowski spacetime, we obtain that complete spacelike immersed surfaces are graphs 
of $1$-contracting maps from the Euclidean space into $\mathbb R$ (Mess' Lemma $1$). In anti-de~Sitter space $\mbox{AdS} \approx \mbox{PSL}(2,{\mathbb R})$, the right action of $\mbox{SO}(2)$ defines a timelike isometric flow
$\Phi^t$ for which the orbit space $(Q^\Phi, g_0)$ is isometric to the hyperbolic plane ${\mathbb H}^{2}$ (in order to be coherent with the previous presentation, we should consider the universal covering $\widetilde{\mbox{SO}}(2)$ acting on $\widetilde{\mbox{SL}}(2,{\mathbb R})$). Corollary~\ref{icon} applies: see Mess' Lemma $6$.

\begin{rac}
\em
In de~Sitter space, there is no everywhere-timelike Killing vector field and so we cannot use Lemma~\ref{increase}. Indeed there are non-injective isometric immersions of complete Riemannian hypersurfaces in de~Sitter space; see \S\ref{S:deS}.
\end{rac}

\begin{rac}
\em
These ideas have been extended to the situation in which $Q^\Phi$ is not
necessarily the orbit space of a timelike Killing vector field, but only
the leaf space
of a $1$-dimensional foliation with timelike leaves. See \cite{HL01}.
\end{rac}

\section{Standard spacetimes}
\label{standard} 

The first main result of this section is Proposition 4,
which shows that the developing image of a spacelike slice
is an embedding and that the time coordinate is proper.
In addition, this section contains the important definition of
\textit{standard spacetimes}, and it is
shown in Proposition 6 that a spacelike slice in an
arbitrary flat spacetime admits a neighborhood that embeds
in a standard spacetime.

The second half of this section handles, more or less completely,
the genus 1 case: a neighborhood of a closed spacelike torus in
a flat spacetime embeds in a complete spacetime or in a ``standard''
example (Definition 2; these are {\em torus spacetimes}
in \cite{Ma84}, \cite{Ca98}). 

\section{Domains of dependence and geodesic laminations}
\label{dod}

This section is the heart of the paper and contains the main theorem
in the flat case, which shows that domains of dependence
(genus at least two) are in one-one correspondence with measured
geodesic laminations.

A key idea introduced in this section is the notion of a 
{\em domain of dependence}.  Mess provides a definition peculiar to  the flat case (Definitions $3$ and $4$).
He also mentions in the introduction a more general definition, characterizing domains of dependence as spacetimes which are maximal with respect to the property that there is a closed (compact without boundary) spacelike hypersurface through each point.
It appeared to us pertinent for the reader's convenience to clarify this
notion which is classical in General Relativity; this is
done below in \S\ref{s.defgh} and \S\ref{s.defmax}.
Domains of dependence are
\emph{maximal globally hyperbolic spacetimes} (Definitions N.5.5 and N.5.11),
while the definition of Mess
restricts to the \emph{Cauchy compact} case (see the final remark of \S\ref{s.defgh}).

The section opens with a discussion of some basic
causality notions;
domains of dependence are defined and
Proposition 11 describes the structure of the causal horizon.
(Note that several of the conclusions of Proposition 11 hold
in greater generality; see for example \cite[Ch. 6]{HE73}). 
Proposition 12 is one half of the main theorem: given a hyperbolic surface
and measured geodesic lamination, there is a corresponding flat spacetime.
This result is a Lorentzian version of the {\em grafting} operation for
complex projective surfaces, and readers unfamiliar with this construction
may profit from first reading the details in \cite{KT92}.
Mess' proof is long and contains a lot of implicit information on the structure
of these examples.  
Proposition 13 is the other half: an arbitrary domain of dependence
(in fact, its causal horizon) determines a measured geodesic lamination,
inverse to the correspondence in Proposition 12.
This is followed by a short interlude discussing the
identification of the moduli space of domains of dependence with
the (co)tangent bundle of Teichm\"uller space, the symplectic structure on it,
and the Chern-Simons reduction due to Witten. 
The section concludes with Propositions 14 and 15 which show that the group
action on the causal horizon is dynamically complicated; Mess uses this
to characterize
when domains of dependence embed in larger spacetimes.  We recast this
result in terms of the Cosmic Censorship Conjecture in \S\ref{SS:CC} below.

The following theorem summarizes the results of \S5 that lead to
the parameterization in terms of measured geodesic laminations.

\begin{thm}
\label{flatdependence}
Let $M$ be a flat oriented Cauchy compact domain of dependence of dimension $2+1$.
Up to a reversal of time orientation, $M$ is future complete. Moreover, $M$ is the quotient of
a convex domain $D$ by a discrete group of isometries $\Gamma$. Let $g$ be the genus of the Cauchy surfaces. Then, $g$ is positive. If $g>1$, the linear part of the holonomy group $\Gamma$ is Fuchsian. If $g=1$, $M$ is either complete, or a Torus Universe.
\end{thm}

For the definition of Torus Universes, see \cite{Ca98}. In the preprint, Mess calls them ``standard spacetimes'' (with toral Cauchy surfaces). They are called {\em Misner spacetimes} in \cite{Ba05}. The theorem above is generalized in \cite{Ba05} to higher dimensions, and by replacing the hypothesis ``Cauchy compact'' by the requirement that \emph{one\/} Cauchy surface is complete.

\subsection{Global hyperbolicity}
\label{s.defgh}

The notion of \textit{domain of dependence} is related to the causal notion
of global hyperbolicity. This is essentially the strongest assumption one can
make regarding the causal properties of a given spacetime; an intermediate
but fundamental notion, necessary for the definition of global hyperbolicity,
is \textit{strong causality}.  We remind the reader of the
necessary definitions.

A causal (resp. timelike) curve is an immersion $c: I \subset {\mathbb R} \rightarrow M$ 
such that for every $t$ in $I$ the derivative $c'(t)$ is causal (resp. timelike). 
This notion extends naturally to non-differentiable curves (see \cite{BEE96}).
Such a curve is \emph{extendable\/} if there is another causal curve 
$\hat{c}: J \rightarrow M$ and a homeomorphism $\varphi: I \rightarrow K \subset J$ 
such that $K \neq J$ and $c$ coincides with $\hat{c} \circ \varphi$. The causal curve 
$c$ is \emph{inextendable\/} if it is not extendable.


\begin{rac}
\em
Spacelike hypersurfaces are locally strictly achronal; non-timelike hypersurfaces are locally achronal.
\end{rac}

\begin{definition}
Let $x$ and $y$ be two points in $M$, with $y$ in the future of $x$. The common past-future region $U(x,y)$ is the intersection between the past of $y$ and the future of $x$.
\end{definition}

The domains $U(x,y)$ form the basis for a topology on $M$, the so-called \emph{Alexandrov topology\/} (see \cite{BEE96}). Observe that every $U(x,y)$ is open
for the manifold topology. The converse in general is false:

\begin{definition}
If the Alexandrov topology coincides with the manifold topology, $M$ is said to be {\em strongly causal}.
\end{definition}

\begin{rac}
\em
If $M$ is strongly causal, every open domain $U \subset M$ equipped with the 
restriction of the ambient Lorentzian metric is strongly causal.
\end{rac}

\begin{propo}[Proposition $3.11$ of \cite{BEE96}]
\label{pro.strongconvex}
A Lorentzian manifold $M$ is strongly causal if and only if it satisfies the
following property: for every point $x$ in $M$, every neighborhood of $x$ 
contains an open neighborhood $U$ (for the usual manifold topology) which 
is \emph{causally convex,\/} i.e., such that any timelike curve in $M$ 
joining two points in $U$ is actually contained in $U$.\fin
\end{propo}

\begin{definition}[\cite{BEE96}, p. 11]
\label{def.gh}
$M$ is {\em globally hyperbolic} if:

-- it is strongly causal,

-- for any $x$, $y$ in $M$, the intersection between the causal future of $x$ and the causal past of $y$ is compact.
\end{definition}

From now on, we assume that $M$ is strongly causal.

The notion of global hyperbolicity is closely related to the notion of Cauchy 
surfaces: let $S$ be a spacelike surface embedded in $M$.

\begin{definition}
The {\em past development} $P(S)$ (resp. the {\em future development} $F(S)$) is the set of points $x$ in $M$ such that every inextendable causal path containing $x$ meets $S$ in its future (resp. in its past). The {\em Cauchy development} ${\mathcal C}(S)$ is the union $P(S) \cup F({S})$.
\end{definition}

\begin{definition}
If $S$ is acausal and ${\mathcal C}(S)$ is all of $M$, then $S$ is said to be a {\em Cauchy surface}.
\end{definition}

\begin{thm}[\cite{Ge70}]
\label{gerochthm1}
A Lorentzian manifold $M$ is globally hyperbolic if and only if it admits a Cauchy surface. \fin
\end{thm}

\begin{thm}[\cite{Ge70}, Proposition 6.6.8 of \cite{HE73}]
If $M$ is globally hyperbolic, and $S$ is a Cauchy surface of $M$, there is a diffeomorphism
$f: M \rightarrow S \times {\mathbb R}$ such that every $f^{-1}(S \times \{ \ast \})$ is a Cauchy surface in $M$.\fin
\end{thm}

\begin{rac}
\em
There has been some imprecision in the literature concerning the proof of
the smoothness of the splitting of globally hyperbolic spacetimes. 
See \cite{BS03a, BS04, BS05b} for a survey of this question, and a complete proof of the smoothness of the splitting $M \approx  S \times {\mathbb R}$.
\end{rac}

\begin{rac}
\em
A globally hyperbolic spacetime is said to be \emph{Cauchy compact} 
(or \emph{spatially compact\/}) if it admits a closed Cauchy hypersurface. All the Cauchy hypersurfaces are then closed, and every closed embedded spacelike hypersurface is a Cauchy hypersurface. The (well-known) proofs of these assertions are as follows: if $S' \subset M \approx S \times {\mathbb R}$ is a (connected) closed spacelike hypersurface, the projection on the first factor $S$ induces a local homeomorphism $p: S'\rightarrow S$. The compactness implies that $p$ is a covering map. The chronological orientation induces an ordering on the fibers of $p$ which is preserved by lifting of curves. It follows that this covering is trivial, i.e., a homeomorphism. Hence, $S'$ is the graph of some
function $f: S \rightarrow \mathbb R$. 
Let $t: M \rightarrow \mathbb R$ the projection on the second factor. Since $S$ is compact,
$t \circ f$ has bounded image. Since $t$ is strictly increasing from $-\infty$ to $+\infty$ along every inextendable causal curve, it follows that $S'$ is a Cauchy hypersurface.
\end{rac}

\subsection{Maximal globally hyperbolic spacetimes}
\label{s.defmax}

\begin{definition}
An isometric embedding $f: M \rightarrow N$ is a {\em Cauchy embedding} if the image under $f$ of any Cauchy surface in $M$ is a Cauchy surface of $N$.
\end{definition}

\begin{definition}[\cite{HE73}, pp. 251--252]
\label{def.maxigh}
A globally hyperbolic spacetime $M$ is {\em maximal} if every Cauchy embedding $f: M \rightarrow N$ in a spacetime with constant curvature is surjective.
\end{definition}

\begin{thm}[see Choquet-Bruhat-Geroch \cite{CBG69}] \label{T:CBG}
Let $M$ be a globally hyperbolic spacetime with constant curvature. Then there is a Cauchy embedding $f: M \rightarrow N$ in a maximal globally hyperbolic spacetime $N$ with constant curvature. Moreover, this maximal globally hyperbolic extension is unique up to right composition by an isometry.\fin
\end{thm}

\begin{rac}
\em
Actually, this theorem admits a much more general extension: instead of restricting 
to spacetimes with constant curvature, we need only require that the spacetimes under consideration satisfy the Einstein equation (in the void, with 
cosmological constant).
\end{rac}

\begin{rac}
\em
\label{cauchyproblem}
The notion of  global hyperbolicity is linked with  the most usual way
to find  solutions of the  Einstein equation: i.e., solving the associated
Cauchy  problem.   This  approach consists of
considering a manifold $\Sigma$  endowed with a Riemannian metric $\bar{g}$ and
a symmetric $2$-tensor $II$, and trying to find a Lorentzian metric $g$
on  $M=\Sigma\times (-1,+1)$,  such  that $g$  satisfies the  Einstein
equation,  such   that  $\bar  g$   is  the  restriction  of   $g$  on
$\Sigma=\Sigma\times\{0\}$  and such that  $II$ represents  the second
fundamental  form of  $\Sigma=\Sigma\times\{0\}$  in $M=\Sigma  \times
(-1, +1)$.  For the problem to admit a solution, the initial data
$(\Sigma,\bar{g},II)$ must satisfy the \emph{constraint 
equations}  (for geometers,  the contracted Gauss-Codazzi  equations).  Conversely,
the Choquet-Bruhat theorem (\cite{Fo52})  states that every set of initial data
satisfying the constraint  equations leads to a solution, which, by the
nature of the process, is globally hyperbolic.  Furthermore, according to the
Choquet-Bruhat-Geroch Theorem (\ref{T:CBG}) above, there is a  unique maximal 
globally hyperbolic solution (up to isometry). 

\end{rac}

\subsection{Cosmic Censorship}
\label{SS:CC}

Let's reproduce here part of the content of Proposition $15$:

\begin{propo}
\label{flatcc}
Let $M$ be a flat domain of dependence of dimension $2+1$. Assume that the holonomy of $M$ fixes a point
in Minkowski space, or that the associated measured geodesic lamination has no isolated 
leaf. Then, $M$ is a maximal flat spacetime, i.e., any isometric embedding 
$f: M \rightarrow N$ in a flat spacetime $N$ (even not globally hyperbolic)
is surjective.
\end{propo}

The proof of this proposition rests on the study of the dynamical properties of
the group action on the boundary of $D$, where $D$ is the convex domain in
Minkowski space such that $M = \Gamma\backslash{D}$. Actually, nearly
the same proof shows that this proposition holds also for anti-de~Sitter or de~Sitter domains 
of dependence (see \cite[Ch.10]{Sc96}). 

This proposition has a definite physical flavour: indeed, it resolves a special case of the Strong Cosmic Censorship Conjecture:

{\bf Strong Cosmic Censorship Conjecture} \cite{Pe82}:
\begin{em}  
Let $M$ be a compact manifold. Then, for every set of generic Cauchy data $(M, \bar{g})$,
the maximal globally hyperbolic domain of dependence defined by these data is a maximal spacetime in the $C^2$ category.
\end{em}

Since measured geodesic laminations without isolated leaves are dense among measured 
geodesic laminations, Proposition~\ref{flatcc} is essentially the proof of 
Cosmic Censorship in the flat $2+1$-dimensional case
(modulo the verification that the map associating a measured
geodesic lamination to Cauchy data is open -- this fact does
not appear in Mess' paper, nor has a proof been given elsewhere in the
literature to the best of our knowledge).
The hypothesis ``without isolated leaves'' is essential,
since, as explained by Mess,
it is easy to extend 
the domain of dependence outside $D$
if the measured geodesic lamination contains an
isolated leaf. This is one way to see that the Cosmic Censorship Conjecture
fails if the ``generic'' requirement is dropped.

For more details, see \S2.1 of \cite{An04}.

\begin{rac}
\em
The content of the paper \cite{IKSK96} is the proof of Proposition~\ref{flatcc}
for the flat case in any dimension, but only in the case where the holonomy 
preserves a point in Minkowski space. This case follows actually quite easily from
the minimality of strong stable foliations of geodesic flows of negatively curved
Riemannian manifolds: this is precisely what is re-proved in \cite{IKSK96}.
\end{rac}

\subsection{CMC time}
In the introduction, Mess indicates that he has nothing to say about
``the interesting question of foliating the manifolds we consider by surfaces of constant mean curvature [...]''.
At the end of \S5, he mentions this problem once again,
cites the work of Moncrief \cite{Mo89},
and leaves the problem of existence and uniqueness of such a foliation
as a question.
This question has been answered recently; see \cite{BBZ03}, \cite{BBZ05}, and \cite{BZ04} for a sketchy proof in the de~Sitter case. 
We summarize these results here, postponing a more 
detailed discussion of CMC foliations and their asymptotics
(focusing on the flat case) to \S\ref{S:CMC}.

\begin{definition}
A function $t: M \rightarrow {\mathbb R}$ on a domain of dependence $M$ is a {\em CMC time} if 
$t$ is a time function, i.e., strictly increasing along every causal curve, and every fiber $t^{-1}(s)$ is a spacelike hypersurface with constant mean curvature $s$.
\end{definition}

\begin{rac}
\em
The sign of the mean curvature is important. The convention for defining the second fundamental form of $S_s = t^{-1}(s)$ is to take as normal vector field along $S_s$ the future-oriented one. 
It follows from the maximum principle that, for a given spacetime $M$,
if a CMC time function exists, then it is unique, and is as regular as the Lorentzian manifold $M$.
See \cite{BBZ05} for more details.
\end{rac}

\begin{thm}[Theorem $12.1$ in \cite{Ba05}; see also \cite{An02a}]
Every flat future-complete Cauchy-compact domain of dependence $M$ admits 
a CMC time function $t: M \rightarrow (-\infty, 0)$, except if $M$ is complete.
\end{thm}

The de~Sitter case is also valid (see \cite{Sc99} for the definition of parabolic de~Sitter domains of dependence).

\begin{thm}[Theorem $5.4$ in \cite{BZ04}]
Every future-complete Cauchy-compact domain of dependence of constant positive curvature and dimension $2+1$ admits a CMC time function $t: M \rightarrow (-\infty, -2)$, except if $M$ is complete or parabolic.
\end{thm}

The anti-de~Sitter case is the easiest to formulate:

\begin{thm}[\cite{BBZ05}]
Every Cauchy-compact domain of dependence of constant negative curvature and dimension $2+1$ admits a CMC time function $t: M \rightarrow (-\infty, +\infty)$.
\end{thm}

\section{The case of de~Sitter space} \label{S:deS}

In \S6, Mess constructs standard de~Sitter spacetimes $S \times \mathbb{R}$
in terms of (complex) projective structures on $S$ (Proposition 16).
The corresponding classification theorem (that every 
de~Sitter domain of dependence comes from this construction) 
is only proved in special cases: $1+1$-dimensional de~Sitter spacetimes
in Proposition 17 and the genus one $2+1$-dimensional case in
Proposition 18.

There are two standard parameterizations of the space $P(S)$ of
projective structures on $S$; the first is the by the bundle
of holomorphic quadratic differentials over Teichm\"uller space,
where a quadratic differential comes from the Schwarzian derivative
of the developing map of a projective structure.   Mess introduces
these notions briefly, and shows how the Bers embedding gives
a holomorphic section of the bundle $P(S) \to \xteich(S)$.
His assertion that the holonomy map from $P(S)$ to the $PSL(2,\mathbb{C})$
representation variety is a holomorphic submersion is due
to Hejhal \cite{He75}.

The second parameterization was given by Thurston and is
the one used by Mess in the construction of the standard de~Sitter examples.
The key element in the Thurston parameterization is
the {\em grafting} operation, which produces a new projective structure
on $S$ from a given hyperbolic
structure and measured lamination.  In the simplest case that
the measured lamination is a simple closed curve $\gamma$ with weight $t > 0$,
grafting amounts to inserting a projective annulus along $\gamma$
of height $t$.   This extends by continuity to any measured lamination,
defining a map
$$
	\Theta: \xteich(S) \times \xml \to P(S),
$$
which turns out to be a homeomorphism (an exposition of Thurston's
proof is given in \cite{KT92}).

One way to obtain the inverse $\Theta^{-1}$ 
is by ``thickening'' the developing map
$d : \tilde{S} \to \xcp^1$
of a projective structure to an immersion
$D : \tilde{S} \times (0,\infty) \to \mathbb{H}^3$
equivariant with respect to the holonomy 
in $PSL(2,\mathbb{C}) \cong {\rm Isom}^+(\mathbb{H}^3)$.
The measured geodesic lamination is then obtained as a bending
lamination on the frontier of the image of $D$,
where some care is required if $d$ is not an embedding.

After some introductory material on the Klein model of
de~Sitter space and projective structures, Mess gives
the construction of $D$ in terms of so-called
{\em maximal balls} in the universal cover $\tilde{S}$.
There are some minor inaccuracies in this discussion,
such as the assertion that the curves foliating overlapping maximal balls 
develop to circular arcs in $\xcp^1$ with the same endpoints.
We recommend Kulkarni and Pinkall's later paper \cite{KP94}
as the best place to read the details of the maximal
ball construction (the cited paper \cite{KP86} is quite 
different).   Also note that the reference to \cite{Ap81}
was probably intended to have been one of Apanasov's later papers,
e.g. \cite{Ap88}.

Then, the projective dual of the map $D$ is used to 
produce an equivariant immersion $q : \tilde{S} \times (0,\infty) \to \mathbb{S}^3_1$.
These are the {\em standard de~Sitter spacetimes}, and their basic
properties are given in Proposition 16.
Mess conjectures that every $2+1$-dimensional de~Sitter spacetime 
which is a small neighborhood of a closed oriented spacelike surface
embeds in a unique standard spacetime.
This was shown by Scannell in his thesis \cite{Sc96},\cite{Sc99}.

It is well-known that the developing map of a projective structure
on a surface need not be an embedding or even a covering of its
image \cite{Ma69}, and so the same is true for the developing maps
of standard de~Sitter spacetimes.   This is in contrast to the
flat and anti-de~Sitter cases (as discussed in \S\ref{realization} above),
and represents the primary difficulty in
extending the classification theorems to the de~Sitter case.
This is overcome by defining domains of dependence, etc.,
solely in terms of the causal structure induced by the 
developing map, and not in terms of the geometry of de~Sitter space itself.

The maximal ball construction for projective surfaces works in any
dimension when a manifold is equipped with a flat conformal
(M\"obius) structure \cite{KP94}.   It turns out that the classification
theorem extends as well:

\begin{thm} \cite{Sc99}
Let $M$ be a closed $n$-manifold.
Then there is a bijective correspondence between the moduli
space of flat conformal structures on $M$ and the moduli space
of maximal $n+1$-dimensional de~Sitter domains of dependence
homeomorphic to $M \times \mathbb{R}$.
\end{thm}

The second half of \S6 in Mess' paper contains an
approach to his conjecture.   The outline he proposes
is completely different than the proof eventually given by Scannell,
and it would be interesting to know if it can be made
to work, despite being special to the $2+1$-dimensional case.
More precisely, if one has a (time-oriented) $2+1$-de~Sitter spacetime 
which is a neighborhood of a closed oriented spacelike surface $S$,
then there is an associated holonomy representation
$hol : \pi_1(S) \to PSL(2,\mathbb{C})$.  To show the spacetime
comes from a projective structure on $S$ via the construction
above, it would help to know that $hol$ in fact coincides with
the holonomy of some projective structure.  For this, Mess
appeals to the following famous ``prescribed monodromy'' result:

\begin{thm} \label{T:GKM}
Let $S$ be a closed oriented surface of genus at least two,
and let $h: \pi_1(S) \to PSL(2,\mathbb{C})$ be a homomorphism.
Then $h$ is the holonomy representation of some complex projective
structure on $S$ if and only if (1) $h$ lifts to $SL(2,\mathbb{C})$
and (2) $h$ is non-elementary.
\end{thm}

Mess cites Gallo's announcement \cite{Ga89} for this, though the
complete proof was not made available until several years
later, in joint work with Marden and M. Kapovich \cite{GKM00}.
Gallo's proposed proof of Theorem \ref{T:GKM} is based on the
existence of a pants decomposition of $S$ with the property that
the holonomy of each pair of pants is quasi-Fuchsian.
Mess argues from this that (after a small change in the holonomy)
one can assume that the holonomy of each pair of pants is Fuchsian.
But even with this assumption it is not clear how to produce
a projective structure on $S$ without appealing to the more 
general convexity results in \cite{Sc99}.
There is also an issue with changing the holonomy,
as is pointed out in the proof of Proposition 18.

In any case, Mess succeeds in proving the conjecture
for $1+1$-dimensional spacetimes (Proposition 17)
and genus one $2+1$-dimensional spacetimes (Proposition 18).

\begin{rac}
\em
Tan's UCLA thesis \cite{Ta88} and a (never published) preprint
of Gallo, Goldman, and Porter \cite{GGP} are cited for the classification
of projective structures with holonomy in $PGL(2,\mathbb{R})$;
for the former, see \cite{Ta94a} and for the latter \cite{GKM00}
(cf. \cite{Go87}).
\end{rac}

\section{Anti-de~Sitter manifolds}
\label{sec.ads}

In \S7, Mess proves a classification theorem for anti-de~Sitter
spacetimes that is analogous to the classification in the flat case;
we summarize this as follows:

\begin{thm}
\label{adsdependence}
Let $M$ be an oriented Cauchy compact domain of dependence locally modeled on anti-de~Sitter space. Then, $M$ is the quotient of
a convex domain $H$ of $\mbox{AdS}$ by a discrete group of isometries $\Gamma$. The genus $g$ of the Cauchy surfaces is positive, and if  $g>1$, the holonomy morphism $\rho = (\rho_L, \rho_R) \rightarrow \mbox{PSL}(2,{\mathbb R}) \times \mbox{PSL}(2,\mathbb{R})$ is a pair of Fuchsian representations.
All pairs $(\rho_L, \rho_R)$ of 
Fuchsian representations are realized.
\end{thm}

In the second half of the section
(immediately following Proposition 21) Mess recasts these results in the
language of {\em earthquakes}, observing that they give an alternative
proof of Thurston's theorem that any two points of Teichm\"uller
space differ by a left earthquake.

The section opens with a discussion of the projective model of
anti-de~Sitter space and identifies the identity component of its
isometry group as
$PSL(2,\mathbb{R}) \times PSL(2,\mathbb{R})$.
Proposition 19 shows that the holonomy of a closed spacelike surface
(implicitly with genus $g>1$)
in an anti-de~Sitter spacetime maps to a Fuchsian group in each factor.
(The proof is analogous to Proposition 1 for the flat case, by arguing
that the Euler classes must be $2-2g$ and then appealing to Goldman's
Theorem).   Thus the holonomy determines a pair
$(x,y)$ in $\xteich(S) \times \xteich(S)$, and the left and right
representations are conjugated by some homeomorphism
$h:\xrp^1 \to \xrp^1$.
The graph of $h$ is a closed achronal topological circle at infinity.
The convex hull $X(\phi)$ of this graph quotiented by the holonomy is an
anti-de~Sitter spacetime $X(x,y)$.


The existence of the spacetime $X(x,y)$ is
the main part of Proposition 20.   It is also claimed
that $X(x,y)$ is uniquely determined by $x$, $y$, and the
fact that its boundary is spacelike, locally convex and without extreme points.
No proof is given however, with the discussion of Proposition 20
ending abruptly with the statement
``Now for the uniqueness of $X(\phi)$'' before Lemma 6.

Let $S$ be a closed spacelike surface in a spacetime $M$
locally modeled on $\mbox{AdS}$.
Then, the inclusion $S \subset M$ is incompressible, and the restriction
to $\widetilde{S}$ (a lift in $\widetilde{M}$ of $S$) of $\mathcal D$ is
injective (see Lemma 6 and compare Lemma 1 from the flat case).
These statements follow from the arguments given above in \S\ref{realization}.
But something more is needed (given by Mess as Lemma $7$):
the restriction of $\mathcal D$ to $\widetilde{S} \approx {\mathbb D}^2$ extends to a continuous map $\overline{\mathcal D}$ from the closed disc $\overline{\mathbb D}^2$ into the closure of $\mbox{AdS} \approx X$ in $\xrp^3$, such that the image by $\overline{\mathcal D}$ of the boundary $\partial{\mathbb D}^2$ is a topological circle in the hyperboloid $Q$. 
Mess proposes two sophisticated proofs of this statement.
Actually, there is a very simple proof of this fact,
using the natural conformal embedding of $\mbox{AdS}$ into
${\mathbb S}^2 \times {\mathbb S}^1$: see \cite{BBZ05} or \cite{Ba05a}.

During the proof of Lemma $7$, Mess shows that the Cauchy development $H(A)$ of $A$ only depends on $A_\infty$: write it as $H(A_\infty)$.
This leads to Proposition 21: any closed spacelike surface $S$ in an anti-de~Sitter
spacetime admits a neighborhood which can be embedded in a
domain of dependence $\Gamma\backslash{H}(A_\infty)$.

Actually, Proposition $21$ starts with a claim that the genus of such a surface has to be greater than one. This claim is false:
there are indeed anti-de~Sitter domains of dependence with toral Cauchy surfaces: the so-called Torus universes (which seem to have been first defined in \cite{Ma84}; see also \cite{Ca98}, \cite{Ba05a}, and the last section of \cite{BBZ05}).
Mess proposed two ``proofs'' for this statement. The first: ``We could argue that the holonomy cover $\hat{F}$ is quasi-isometric to a hyperbolic plane, and so the ball of radius $r$ in $\hat{F}$ grows exponentially with $r$, while the cover of a torus has only polynomial volume growth''. But he doesn't argue more. ``Alternatively, the left holonomy defines a homomorphism from ${\mathbb Z} \oplus {\mathbb Z}$ into an abelian subgroup [...] (with, after perturbation) cyclic image. The right holonomy is topologically conjugate to the left holonomy so the holonomy has kernel which contradicts lemma $6$.''.  The error is that left and right holonomies are not conjugate, only semi-conjugate; the proof of the conjugacy involved Goldman's theorem proving the representations into $\mbox{PSL}(2,{\mathbb R})$ with maximal Euler
  number are all Fuchsian, and hence is not valid for $g=1$. Actually, when $g=1$, the topological circle $A_\infty$ is a union of four lightlike segments. It can be interpreted as the graph of a semi-conjugacy, shrinking two intervals of $\xrp^1$ to two points, and expanding the two common extremities of these intervals to intervals. 

As noted earlier, the second half of \S7 reinterprets the geometric data on
the boundary of the convex hull $X(\phi)$ as earthquakes, reproving Thurston's
Earthquake Theorem.
Then, Mess considers the case where $S$ is not compact.
He considers the example of an earthquake ``taking a complete hyperbolic surface with cyclic fundamental group and one cusp to one of the components of the complement of a closed geodesic in a hyperbolic surface with cyclic fundamental group and no cusp''. This example corresponds to \emph{extreme black-holes\/} as described in \cite{Ba05b}. He also considers the example of the thrice punctured sphere, which is studied in more detail in \cite{BB05}.

The section closes with a number of interesting questions:
the volumes of the domain of dependence and of the convex hull define two maps
on $\xteich(S) \times \xteich(S)$.
``How do they behave? Are they related, perhaps asymptotically, to such invariants of a quasi-Fuchsian groups as the volume of the convex hull and the Hausdorff dimension of the limit set?''.
It is worth noting that the meaning of this question can be clarified thanks to Benedetti-Bonsante Wick rotation (see the penultimate remark in \S\ref{S:wick}
and also \S4.8.2 of
\cite{BB05} for specific computations in this direction).

Mess then asks several questions about the geometry of the boundary of the
convex core in both the anti-de~Sitter and hyperbolic settings.
Concerning
the problem of prescribing the induced metric on the boundary of the 
convex core, in the hyperbolic setting, he mentions that existence 
follows from work of Epstein and Marden.  This existence statement
is also a direct consequence of a result of Labourie \cite{La92a}.

As for prescribing the bending lamination, still in the hyperbolic setting, 
the existence has been proved recently by Bonahon and Otal
\cite{BO04} for convex co-compact manifolds with incompressible boundary
(and, for general convex co-compact
manifolds, by Lecuire \cite{Le02}). However, for prescribing either the
induced metric or the bending lamination, the uniqueness remains unknown,
basically because we do not know whether infinitesimal rigidity holds,
i.e. whether any first-order deformation of a quasi-Fuchsian manifold
induces a non-zero variation of the induced metric (resp. bending
lamination) on the boundary of the convex core. However, in the hyperbolic
setting, Bonahon \cite{Bo96a} proved that the two infinitesimal
rigidity questions, concerning the induced metric and concerning the bending
lamination, are equivalent. Bonahon \cite{Bo05a} also gave
a careful analysis of what happens near the ``Fuchsian locus'', showing 
that uniqueness does hold there.  Series \cite{Se04} proves
uniqueness in the case of once-punctured tori.

Mess also asks whether it is possible to prescribe the induced metric
on one boundary component, and the bending lamination on the other. 
Lecuire \cite{Le06} has recently proved the existence part of this
statement (again in the hyperbolic case). In the same preprint
he also gives a positive answer to the existence part of another
of the questions asked by Mess: whether it is possible to prescribe
the conformal structure on the upper surface at infinity and the
measured bending lamination on the lower boundary of the convex core;
the uniqueness holds when the prescribed bending lamination is
supported on simple closed curves.
The first result is largely a consequence of the second, which is
proved using recent results of Bromberg \cite{Br04} on the 
deformations of complete hyperbolic metrics with cone singularities
along closed curves.

Mess also asks whether a quasifuchsian hyperbolic manifold is 
uniquely determined by the conformal structure on the upper surface
at infinity along with the induced metric on the lower boundary
of the convex core. Here again, the existence should hold, and 
should follow directly from the argument that he mentions above,
using the fact that the induced metric on the lower boundary of
the convex core is quasi-conformal to the conformal structure
on the lower surface at infinity (according to \cite{EM87}).

Another known result in this same vein is that 
one can prescribe a conformal structure at infinity and the 
bending lamination associated to that conformal structure
and these data determine the quasi-Fuchsian group uniquely; see \cite{SW02}.

The article goes on to ask similar questions in the anti-de Sitter
setting. The questions on prescribing the induced metric and/or 
the measured bending lamination on the boundary of the convex 
core remain open. However the arguments given in section 7 give
an answer to one of the questions: whether a maximal globally
hyperbolic AdS manifold is
uniquely determined by its left holonomy (say $h_l$) and the 
induced metric on the upper boundary of its convex core.
By Thurston's earthquake theorem, there is a unique measured
lamination $\lambda_+$ such that $h_l$ is the image of $\mu_+$
by the left earthquake along $\lambda_+$. The right holonomy 
$h_r$ is then the image of $\mu_+$ by the right earthquake along
$\lambda_+$. A similar argument also shows that a maximal globally
hyperbolic AdS manifold is uniquely determined by the left holonomy and the
measured bending lamination on the upper boundary of the convex 
core.

It is perhaps worth remarking that one can replace in these questions the 
convex core by a slightly larger domain, with smooth and strictly convex
boundary. In the  hyperbolic setting, the induced metric now has curvature
$K>-1$, and it is actually possible to obtain any pair of such metrics
uniquely (see \cite{Sc06}). For the same convex domains the smooth analogue
of the bending lamination is the third fundamental form, which is a smooth
metric with $K<1$ and contractible closed geodesics of length larger than
$2\pi$, and, here again, any pair of such metrics can be uniquely prescribed
(also \cite{Sc06}). 

The anti-de~Sitter case remains more elusive than the quasi-Fuchsian case, though some
of the results known in the hyperbolic setting might extend to anti-de~Sitter
manifolds.  In particular, this is the case for the analysis of
\cite{Bo05a}, and also for \cite{Sc06} (for which however the proof
has to be adapted and the dictionary between quasi-Fuchsian and anti-de~Sitter manifolds
has to be extended a little). 

\begin{rac}
\em
Theorem~\ref{adsdependence} has been extended recently (in \cite{BB05} and \cite{Ba05a}) to the case of domains of dependence of dimension $2+1$ admitting complete Cauchy surfaces. But the classification of Cauchy compact domains of dependence in higher dimensions remains an open problem.
\end{rac}

\begin{rac}
\em
Another proof that the moduli space of anti-de~Sitter domains of dependence
is parameterized by two copies of Teichm\"uller space is given 
in \cite{KS05}.  This work is based on
differential-geometric ideas and uses a result of Labourie \cite{La92}
on mappings between hyperbolic surfaces.
\end{rac}

\begin{rac}
\em
A striking fact is that geometric ideas presented in Mess' preprint provide an essential complement to the abundant literature devoted to \emph{BTZ (multi) black-holes\/} (see \cite{BTZ92, BHTZ93, ABHP96, ABBHP98, ABH99, Br96a, Br00, Ca98}, etc... This matter is developed in \cite{Ba05a,Ba05b}.
\end{rac}

\begin{rac}
\em
The analogy between representations of surface groups in the isometry
group of anti-de~Sitter space and quasi-Fuchsian groups has been
pursued and generalized in recent work of Labourie,
in which a definition of a quasi-Fuchsian group in a general Lie group $G$
is given.
Such a notion is characterized by the choice of a specific
$H=SL(2,\mathbb R)$ in $G$.
Surface groups which factor through a cocompact subgroup of $H$ are called
{\em Fuchsian}. 
Then {\em quasi-Fuchsian} (or {\em Anosov}) subgroups are deformations of
these subgroups which satisfy a dynamical condition, and the classical
stability theorem of
dynamical systems ensures that the space of quasi-Fuchsian representations is
open.  
It follows from the definition that a quasi-Fuchsian group generates a
continuous limit curve at a suitable infinity which is a H\"older equivariant
mapping from the boundary at infinity of the surface to a parabolic
quotient.
One may make the link with the situation described by Mess 
following Proposition 19, which gives the curve at infinity
in $\mathbb{RP}^1\times\mathbb{RP}^1$. 

It turns out that two classes of Lie groups have some
special and interesting features in this context: real split groups and isometry groups of
hermitian symmetric spaces of tube type. 
In these cases, the limit curve enjoys an extra {\em positivity} condition (see Fock-Goncharov \cite{FG05}, and Burger-Iozzi-Labourie-Wienhard \cite{BILW05}) which, in the case of $SO(2,2)$ (which belongs to both classes), amounts to the condition that it be a
spacelike curve. 
In both cases, it also turns out that quasi-Fuchsian groups fill out a whole component of the space of representations.
\end{rac}

\begin{rac}
\em
The ideas involved in Mess' proof of the Earthquake Theorem have been
extended recently to two contexts. In \cite{BS06}, the 3-manifolds considered
are AdS manifolds with ``particles'' -- cone singularities along maximal
time-like lines, with angle less than $\pi$ -- but which are globally
hyperbolic. These manifolds happen to have a well-defined convex core,
which leads to a version of the Earthquake Theorem for hyperbolic
surfaces with cone singularities of angle less than $\pi$. 
In \cite{BKS06}, the same idea is used for multi-black
holes, which are not globally hyperbolic. These manifolds also have
a kind of convex core, and considering its geometry leads to an
extension of the Earthquake Theorem for hyperbolic surfaces with
geodesic boundary: the measured laminations on the
{\it interior} of those surfaces act simply transitively by right
earthquakes
on the {\it enhanced} Teichm\"uller space. The enhanced
Teichm\"uller space -- recently introduced by Fock -- contains
an open dense subset which is a $2^n$-cover of the ``usual''
Teichm\"uller space, where $n$ is the number of boundary components.
\end{rac}

\section{Classification of spacetimes} \label{S:classif}

This section starts with a report on results about complete
flat Lorentzian manifolds, in particular the proof (Proposition 24) that the
linear holonomy is either solvable or discrete.  This result
appeared earlier (with a different proof that appeals to
\cite{Au63}) in Fried-Goldman \cite{FG83}.

This is followed by ``Mess' Theorem'' (Proposition 25):
the fundamental group of a closed surface cannot act properly on
all of Minkowski space.  The proof essentially relies on the
``Cosmic Censorship'' machinery from \S5.   New proofs have been
given recently by Goldman-Margulis \cite{GM00}, and Labourie \cite{La01}. 

A consequence of this result is that the fundamental group of
a complete, non-compact, flat Lorentz $3$-manifold is solvable or free.
Mess goes on to say that ``it seems plausible that a complete Lorentz
manifold with free fundamental group is diffeomorphic to the interior of
a (possibly non-orientable) handlebody'', but this question remains open.
On the other hand, the cited conjecture of Margulis that a complete
Lorentz manifold with finitely generated free fundametal group has discrete
and purely hyperbolic linear holonomy has been disproven;
Drumm \cite{Dr92} constructs
affine deformations of any free Fuchsian group which act
properly on Minkowski space.

At the heart of the Goldman-Margulis proof of Mess' Theorem
is the result that the
Margulis signed length (see \cite{GM00}) is the Hamiltonian flow of the length in 
Teichm\"uller space; this is conjectured by Mess 
(and ``left to the reader'') in the extended discussion 
of quantization and Teichm\"uller theory 
following Proposition 25.  Quantization of Teichm\"uller spaces
has become a classical theme in the last several years.  
The initial breakthrough was due to Chekhov and Fock \cite{FC99},
who quantized the algebra of observables given by the shear
coordinates on open surfaces.
By developing this technology, invariants of hyperbolic
$3$-manifolds and more generally of $3$-manifolds 
equipped with $PSL(2,\mathbb{C})$ characters have been
constructed (see Kashaev \cite{Ka98},
followed by works of Baseilhac and Benedetti \cite{BB04}, \cite{BB05a}
and Bonahon and Liu \cite{BL04}).
As yet, however, the relationship between anti-de~Sitter and physical
2+1 gravity is not completely clear.

Mess' suspicion that the spectrum of the Weil-Petersson Laplacian on 
moduli space is bounded away from zero has been confirmed 
by McMullen (in his paper showing moduli space is K\"ahler hyperbolic \cite{Mc00}).

Next Mess proves the following theorem (Proposition $26$):

\begin{thm}
\label{flatcompact}
Let $M$ be a compact flat $2+1$ spacetime with spacelike boundary (maybe empty). Then,
either $M$ is complete (if the boundary is empty), or diffeomorphic to $S \times [0,1]$ where
every $S \times \{ t \}$ is spacelike.
\end{thm}

Actually, Mess' proof implies slightly more: in the non-closed case, the interior of $M$ is globally hyperbolic. This theorem is a natural extension of Carri\`ere's Theorem establishing the completeness of closed flat spacetimes (\cite{Ca89}). 


It seems possible to simplify the proof by using the cosmological time (see \cite{AGH98}): once one has obtained by Koszul's argument that (in Mess' notation) $\hat{M}'$ is the quotient of an open domain, then one can show that if $M$ is not complete, then the cosmological time of $M'$ is regular (up to time reversing); it follows directly from \cite{AGH98} that $M'$ is globally hyperbolic.

This result is followed by an analogous statement concerning time-oriented
locally anti-de~Sitter compact spacetimes with spacelike boundary (Proposition 27).
Observe that the conclusion of the last statement of this proposition
can be greatly simplified by the result of Klingler \cite{Kl96}:
these spacetimes are either globally hyperbolic
(Mess says that ``they embed in a domain of dependence'')
or complete
(since according to \cite{Kl96}, closed anti-de~Sitter spacetimes are complete).
Indeed, the last two pages of the proof concern the closed case and
can therefore be skipped in favor of \cite{Kl96}.
Observe also that the proof has to be corrected, since Mess starts the proof arguing that the boundary components have genus greater than one
according to Proposition $21$, and we have seen that this is not correct.

Theorem~\ref{flatcompact} implies that the interior of a compact flat $2+1$-spacetime is globally hyperbolic. The analogous statement for anti-de~Sitter spacetimes
holds when the boundary is non-empty, but is false in the closed case:
indeed, quotients of the entire anti-de~Sitter space are never globally hyperbolic.

\begin{rac}
\em
In order to conclude the case with boundary,
Mess argues that the action of the surface group on a connected component
of the horizon of a domain of dependence in $\mbox{AdS}$ is topologically conjugate
to the action on the horizon of a flat domain of dependence to which the
``Cosmic Censorship'' principle applies.
This remark is a precursor to the ``canonical Wick rotation'' discussed in
\S\ref{S:wick}.
\end{rac}

Mess concludes with some comments on the anti-de~Sitter case which are no
longer relevant in light of Klingler's Theorem; later work on the
classification of closed anti-de~Sitter spacetimes can be found in 
\cite{Sa97}, \cite{Sa00}, and \cite{Ze98}.
The questions posed concerning de~Sitter spacetimes have also been
resolved; Klingler's result \cite{Kl96} (see also \cite{Mo96}) implies
that there can be
no closed de~Sitter manifold (any dimension) as predicted by Mess.
The bounded case is a bit more delicate; Scannell's thesis \cite{Sc96} gives
necessary and sufficient conditions for a $2+1$-dimensional de~Sitter domain
of dependence to be maximal.   There are even counterexamples 
to the statement given by Mess in the $1+1$-dimensional de~Sitter case
\cite[Figures 1 and 2]{Sc99}.

At the very end, Mess mentions the higher dimensional (flat) case. He provides an 
argument proving that the holonomy of spacelike hypersurfaces 
is always discrete: this is what we explain here in \S\ref{standard}.
In the flat case, it allows one to avoid the use of Goldman's Theorem: this 
is detailed in \cite{Ba05}.  Indeed, the final claim of the paper
about the flat $3+1$ dimensional case can also be found in \cite{Ba05}.

\section{Canonical Wick Rotation}
\label{S:wick}

\nocite{BB06}

Given an orientable surface $S$ let us denote by $\Mm\Ss_k(S)$ the
space of Lorentzian structures on $S\times\mr$ of constant curvature
$k$ such that $S\times\{0\}$ is a complete Cauchy surface up to the
action of the homotopically trivial diffeomorphisms fixing
$S\times\{0\}$.

When $S$ is compact, the space $\Mm\Ss_k(S)$ has been shown to be
homeomorphic to the cotangent bundle of the Teichm\"uller space of
$S$, provided that the genus of $S$ is at least $2$
(by Mess in the case $k\leq 0$ and Scannell~\cite{Sc99} in the
case $k=1$).
 
In~\cite{BG01}, Benedetti and Guadagnini stressed the role of cosmological
time as a fundamental tool to better understand flat globally hyperbolic
spacetimes classified by Mess. In fact, the cosmological time turns out to
be an important object also in~\cite{Sc99}. A remarkable fact is
that, in both contexts, level surfaces of the cosmological time
are obtained by grafting a hyperbolic surface $F$ (homeomorphic
to $S$) along a measured geodesic lamination $\lambda$.  Moreover,
the Mess parameters are explicitly related to the pair $(F,\lambda)$
(actually they furnish good parameters for the space
$\Mm\Ss_k(S)$).

Similar behavior occurs in the anti-de~Sitter
framework, even if in this case cosmological time is a $\mathrm
C^{1,1}$-function until it reaches the value $\pi/2$. Anyway, also in
this case level surfaces for values $<\pi/2$ are obtained by grafting
a hyperbolic surface along a measured geodesic lamination, and these
data determine the spacetime.\\

When $S$ is not compact, one could try to generalize the parameterizations
of $\Mm\Ss_k(S)$ and $P(S)$.  In fact, the notion of measured
geodesic lamination can be implemented for every hyperbolic surface, and it is
not difficult to see that the Mess and Thurston constructions work as well.
But in this case they do not give rise to a
complete classification of $\Mm\Ss_k(S)$ or $P(S)$ (\emph{i.e.}, there are
globally hyperbolic  spacetimes of constant curvature such that the
cosmological level sets are not obtained by grafting a hyperbolic surface
along a measured geodesic lamination).

There are two natural problems arising from this remark:
\medskip\par\noindent
\emph{1)} To find a more general notion of measured geodesic lamination,
coinciding with the usual one in the compact case, that allows one to
generalize the
Mess and Thurston constructions to obtain complete classifications of
$\Mm\Ss_k(S)$ as well as $P(S)$.
\medskip\par\noindent
\emph{2)} To make explicit the identifications between $P(S)$ and
$\Mm\Ss_k(S)$ for compact $S$ (arising from the Thurston and Mess
parameterizations) in order to see whether they can be
generalized to the non-compact case.
\medskip\par In~\cite{KP94}, Kulkarni and Pinkall introduced the
notion of a measured geodesic lamination on a straight convex set that
allows one to carry out a complete classification of projective structures
on a surface $S$ with non-abelian fundamental group. Actually, they
showed that the Thurston construction could be applied to these more
general laminations and that every projective structure could be
constructed in such a way.\\

In~\cite{BB05}, it is shown that the Mess constructions could be applied
also to these laminations and this leads to a complete classification
of $\Mm\Ss_k(S)$.

In the flat case, the proof is based
on~\cite{Ba05,An02a}, which provide a clear picture of the
universal covering spaces and the linear holonomies of globally
hyperbolic flat spacetimes. On the other hand, in~\cite{Bo06} the
universal covering spaces are classified in terms of measured geodesic
laminations on straight convex sets.

The proofs in the de~Sitter and
the anti-de~Sitter cases are carried over in~\cite{BB05} by
developing the ideas of Mess and Scannell in this more general case,
and by using an explicit map
\[
   \Mm\Ss_0(S)\rightarrow\Mm\Ss_k(S)
\]
that solves question \emph{2)} above.  In fact, such a map is constructed
by developing a \emph{ canonical Wick rotation and rescaling theory}.

Let us briefly introduce these notions.  In general, given a manifold
$M$, a nowhere vanishing vector field $X$, and a pair of positive
functions $\alpha,\beta$, the Wick rotation is an operation
transforming Riemannian metrics on $M$ into Lorentzian metrics that
make $X$ a timelike vector field. Namely, given a Riemannian metric
$g$ the metric $h=W_{(X,\alpha,\beta)}(g)$ obtained by the Wick
rotation of $g$ along $X$ with rescaling functions $\alpha$ and $\beta$
is determined by the following properties:
\medskip\par\noindent
\emph{1.} $X^{\perp_g}=X^{\perp_h}=X^\perp$.
\medskip\par\noindent
\emph{2.} $h|_{X^\perp}=\alpha g|_{X^\perp}$.
\medskip\par\noindent
\emph{3.} $h(X,X)=-\beta g(X,X)$.

\medskip\par Clearly, the Wick rotation can also be regarded as an
operation transforming Lorentzian metrics with $X$ a timelike vector
field into Riemannian metrics.

The rescaling operation similarly depends on a
vector field $X$ and two positive functions $\alpha,\beta$, and acts
on the space of Lorentzian metrics that make $X$ a timelike vector
field. The main difference with respect to the Wick rotation is that
it preserves the signature of the metrics. Namely, the rescaled metric
$h=R_{(X,\alpha,\beta)}(g)$ is determined by properties \emph{1., 2.}
(the same used to define the Wick rotation) and
\medskip\par\noindent
\emph{3'.} $h(X,X)=\beta g(X,X)$.\\

The Wick rotation rescaling theory is developed through the following scheme:

\begin{itemize}
\item Every maximal globally hyperbolic \emph{flat}
spacetime homeomorphic to $S\times\mr$ is proved to be equipped with a $\mathrm
C^{1,1}$ cosmological time $T$ (provided that $\pi_1(S)$ is not abelian).

\item A canonical Wick rotation on
$M(>1):=T^{-1}((1,+\infty))$ directed along the gradient of $T$ is
shown to yield a hyperbolic metric. Moreover this hyperbolic
structure extends to a (complex) projective structure
on the level surface $M(1)=T^{-1}(1)$.

\item A canonical rescaling on $M(<1):=T^{-1}((0,1))$ directed along the gradient
of $T$, and yielding a de~Sitter metric, is pointed out. Such a de~Sitter
spacetime, denoted by $M^{(1)}$, turns out to be maximal
globally hyperbolic, the level surfaces $M(a)$ of $T$ are Cauchy surfaces, 
and its cosmological time is an explicit function of $T$.

\item A canonical rescaling on $M$ directed along
the gradient of $T$ yields an anti-de~Sitter structure denoted
by $M^{(-1)}$. Level surfaces of $T$ are Cauchy surfaces for
$M^{(-1)}$ and its cosmological time is an explicit function of
$T$. It is not maximal but it coincides with the past part of its
maximal extension (that is denoted by $\Nn(\Mm^{(-1)})$).
\end{itemize}

In this context the word \emph{canonical} means that a function
$f:M\rightarrow N$ between two flat globally hyperbolic spacetimes is
an isometry if and only if it is an isometry for the respective Wick
rotated (or rescaled) structures.

The Wick rotation-rescaling theory leads to the following
classification theorem.
\begin{thm}
Let $S$ be a surface with non-abelian fundamental group.
Then the maps
\[
\begin{array}{l}
   \Mm\Ss_0(S)\ni M\mapsto M(1)\in P(S)\\
   \Mm\Ss_0(S)\ni M\mapsto M^{(1)}\in\Mm\Ss_1(S)\\
   \Mm\Ss_0(S)\ni M\mapsto \Nn(M^{(-1)})\in\Mm\Ss_{-1}(S)
\end{array}
\]
are bijective.
\end{thm}
\bigskip\par\noindent

\begin{rac}
\em
In~\S7, Mess relates the classification of
$\Mm\Ss_{-1}(S)$ to the earthquake theory on $S$. In fact,
the holonomy of $M\in\Mm\Ss_{-1}(S)$ is given by a pair of Fuchsian
representations $(\rho_L,\rho_R)$ of $S$ (this makes sense because of
the natural identification of the isometry group of the Klein model of
anti-de~Sitter space with $PSL(2,\mr)\times PSL(2,\mr)$). Moreover, the
universal covering $\tilde M$ of $M\in\Mm\Ss_{-1}(S)$ is a convex
domain in the Klein model of anti-de~Sitter space and its
extension to the boundary (which is canonically identified with
$\xrp^1\times\xrp^1$) is the graph of the unique homeomorphism of
$\xrp^1$ conjugating $\rho_L$ and $\rho_R$.

Let $\lambda$ be the measured geodesic lamination of
$F_L=\mh^2/\rho_L$ such that the left earthquake along it sends $F_L$
to $F_R=\mh^2/\rho_R$. Denote by $F_+$ the surface obtained by a left
earthquake on $F_L$ along $\lambda/2$. Then the level set
$T^{-1}(\pi/2)\subset M$ is obtained by \emph{bending} (in a suitable
sense) $F_+$ into $M$ along the measured geodesic lamination of $F_+$
corresponding to $\lambda$.

When $S$ is not closed, the relationship between anti-de~Sitter geometry and
hyperbolic geometry is more involved. In fact, on the one hand it is
well-known that there are measured geodesic laminations on $S$ that do
not give rise to a genuine earthquake. On the other hand it is no
longer true that the closure of the universal covering of
$M\in\Mm\Ss_{-1}(S)$ is the graph of a homeomorphism (in general it is
just a nowhere-timelike simple curve).

For every measured geodesic lamination on a straight convex set,
a \emph{generalized earthquake} is defined as an injective, but in
general not surjective, map of the straight convex set on which the
lamination is defined into $\mh^2$. The boundary curves of the
universal covering of maximal globally hyperbolic anti-de~Sitter spacetimes can be regarded as the trace
at infinity of such generalized earthquakes. In particular such a curve
is the graph of a homeomorphism if and only if the corresponding
generalized earthquake is surjective.
\end{rac}

\begin{rac}
\em
The class of maximal globally hyperbolic
anti-de~Sitter spacetimes is invariant under the reversal of time-orientation.
On the other hand, the sub-class of those corresponding to
\emph{genuine} measured geodesic laminations on the whole $\mh^2$ is
not invariant under that operation. 
Such a phenomenon is shown by deepening the example 
suggested by Mess of the hyperbolic thrice-punctured sphere 
bent along geodesics joining the punctures. Indeed it is
somehow related to the fact that in general earthquakes
are not surjective.
\end{rac}

\begin{rac}
\em
Let us consider a closed surface $S$ and
$M\in\Mm\Ss_k(S)$.  Since the gradient of the cosmological time $T$ is
a unitary vector field, the following formula holds
\[
   Vol(a,b)=\int_a^b A(t)\mathrm d t
\]
where $V(a,b)$ is the volume of $T^{-1}(a,b)$ and $A(t)$ is the area
of the surface $T^{-1}(t)$.  Since the level set $T^{-1}(t)$ is (up to
rescaling) the grafting of a hyperbolic surface $F$ along a measured
geodesic lamination $g(t)\lambda$ (where $g$ is an explicit positive
function of time) the area of $T^{-1}(t)$ is a function of the
\emph{length} of $\lambda$.

In such a way, formulas are given in~\cite{BB05} which explicitly compute
$V(a,b)$ in terms of the Mess parameters of $M$. In particular for
$k=-1$, the computation allows one to compute the volume of
$T^{-1}(0,\pi/2)$.
\end{rac}

\begin{rac}
\em
Let $S$ be a closed surface of genus at least
$2$.  For any hyperbolic structure $F$ on $S$ and any measured
geodesic lamination $\lambda$ on $F$, denote by
$M^{(k)}(F,\lambda)\in\Mm\Ss_k(S)$ the spacetime corresponding to the
pair $(F,\lambda)$.  For every $t>0$ let us consider the spacetime
$N_t=\frac{1}{t^2}M^{(k)}(F,t\lambda)$ obtained by rescaling the metric
on $M^{(k)}(F,t\lambda)$ by the factor $1/t^2$. By means of the
canonical Wick Rotation rescaling theory one can see that (in a
suitable sense) we have
\[
    N_t\rightarrow M^{(0)}(F,\lambda)
\]
as $t\rightarrow 0$.
\end{rac}

\section{Constant mean curvature foliations}
\label{S:CMC}

\subsection*{Hypersurfaces of constant mean curvature} 
Let $M$ be a spacelike hypersurface with timelike unit normal $T$, 
in a Lorentz spacetime $V$, with metric 
$\la \cdot , \cdot \ra$ and covariant derivative $D$. The second fundamental
form of $M$ is $K(X,Y) = \la T, D_X Y\ra$
and the mean
curvature is $\tr K = \sum K(e_i, e_i)$ where the sum is over an orthonormal frame of
$M$.
$M$ is called a constant mean curvature hypersurface if $d\tr K = 0$. 
Let $\Area(M)$ denote the area of $M$ with respect to the
induced volume element. 
If $\tr K = 0$, $M$ maximizes area with respect to compactly supported
variations, and hypersurfaces with $\tr K = 0$ are therefore called
maximal. A hypersurface with constant mean curvature $\tr K = \tau$ on the
other hand maximizes the action
$$
L(M) = \Area(M) + \tau \Vol(M;M_0)
$$
with respect to compactly supported variations.
Here $\Vol(M,M_0)$ is the volume in $V$ bounded by $M$ and some fixed, suitably
chosen $M_0$.  

The mean curvature of a graph $t = f(x^1, \dots, x^n)$ in $n+1$ dimensional
Minkowski space is given by
$$
\mathcal M [f] = \sum_i D_i \left ( \frac{D_i f}{\sqrt{1 - |Df|^2}} \right ) 
$$
Thus the equation $\mathcal M [f]$ is a quasilinear (non-uniformly) elliptic
operator which therefore satisfies the strong maximum principle. The strong
maximum principle holds with rather mild regularity assumptions; see
\cite{AGH98a}. 

The strong maximum principle allows one to construct barriers for the
variational problem with action $L$, and then results from geometric measure theory
allow one to prove existence of smooth solutions. This was carried out
by Gerhardt \cite{Ge83a} in the setting that is relevant here. 

Let $V$ be a spatially compact maximal globally hyperbolic flat (MGHF)
spacetime.  A point in $V$ is on at most one Cauchy surface with mean
curvature $\tau$, unless $V$ splits as a metric product.
In this case, the spacetime is of
translational type, see below, and the only hypersurfaces
with constant mean curvature are maximal, $\tr K = 0$. 

Let $M$ be a CMC Cauchy surface of dimension $n$ 
in a spatially compact MGHF spacetime $V$. 
By viewing the development $\dev(V)$ of $V$ as a
subset of Minkowski space $\Mink^{n+1}$, the universal covering
$\tM$ of $M$, which is in a natural way a subset of $\dev(V)$, 
can be considered as a CMC hypersurface in Minkowski space
$\Mink^{n+1}$. It follows from the construction that $\tM$ is complete. 
CMC hypersurfaces
in Minkowski space are convex, with non-positive Ricci curvature. Further,
the Gauss map $\phi: \tM \to \HH^n$ is harmonic. This was apparently 
first noticed by
T. K. Milnor \cite{Mi83} in the 2-dimensional case.

\subsection*{CMC foliations of flat spacetimes} 
We consider spatially compact MGHF spacetimes. 
By the classification due to Scannell \cite{Sc01} 
in dimension $3+1$ and Barbot \cite{Ba05} for general dimension, 
these are, up to finite coverings and linear
twisted products (see \cite{Ba05}),
either products of spatially compact MGHF spacetimes with Cauchy
surface of hyperbolic type, with the Euclidean
torus, or of translation type. A spacetime is of translation type if it is
a quotient of a flat spacetime of topology $\Re \times T^n$, with the product
metric. It follows from the above discussion that unless the spacetime is of
translation type, the 
mean curvature foliation defines a time function $\tau$ on $V$. 

Consider spatially compact MGHF spacetimes of dimension $\geq 2+1$ 
with Cauchy surface of hyperbolic type. 
Existence of global CMC foliations was
proved by Andersson, Moncrief and Tromba \cite{AMT97}
in the 2+1 dimensional
case (and any cosmological constant),
assuming existence of one CMC hypersurface.
A proof of existence of global CMC foliations of 2+1 dimensional spacetimes
was given by Barbot et al. \cite{BBZ03} using level sets of the
cosmological time function \cite{AGH98} 
as barriers. A proof of existence
of global CMC foliations of spacetimes of general dimension assuming
hyperbolic spatial topology was given by Andersson
\cite{An02a}. The general case was considered by Barbot
\cite{Ba05}.

In case the spacetime has hyperbolic spatial topology, the CMC hypersurfaces
 are strictly convex, i.e. $K(X,Y) \leq \lambda \la X, Y\ra$ for some
 $\lambda < 0$
and $\Ric < 0$ for the induced metric. 
 In case the spacetime has a torus factor in the sense
discussed above, the constant mean curvature hypersurfaces split as a metric
 product of a Ricci negative factor with a flat
 torus factor. This
follows from the work in \cite{An02a}. 

\subsection*{Asymptotics} 
\subsubsection*{The expanding direction}
Assume the spacetime $V$ 
has hyperbolic spatial topology with Cauchy surface $M$. Let $\Gamma$ be the
linear part of the holonomy representation of the fundamental group of $V$.  
The quotient $\HH^n /\Gamma$ defines a hyperbolic metric $\gamma$ on $M$,
with sectional curvature $-1$. 

Let $g_\tau$ be the metric induced on the 
$\tau$-level set of the CMC foliation of
$V$. As $\tau \nearrow 0$, then $(M, g_\tau)$ expands. However, the scale
invariant metric $\tau^2 g_\tau$ converges in the Gromov sense to 
$$
\lim_{\tau \nearrow 0} \frac{\tau^2}{n^2} g_\tau \to \gamma
$$
More generally, suppose
$V$ has a factor with hyperbolic spatial topology $N$ 
of
dimension $m+1$ and a torus factor of dimension $k$. Let $(N,\gamma)$ be
the hyperbolic geometry determined as above.  In this case, 
the $n=m+k$ dimensional 
scale invariant geometry $(M, \frac{\tau^2}{n^2} g_\tau)$ collapses and 
converges in the Gromov sense,
as $\tau \nearrow 0$, to the $m$-dimensional space $(N,\gamma)$. 

\subsubsection*{The collapsing direction} The development $\dev(V)$ is a
convex subset of Minkowski space. The singularity of $V$ can be identified
with $P = B/\sim$ where $B$ is the boundary of $V$ and $\sim$ is an
equivalence relation; see \cite{Bo05}. It was shown by Benedetti and
Guadagnini \cite{BG01} that the scale invariant 
geometry on the level
sets of the cosmological time function converges in the Gromov sense to that
of $P$, which in the 2+1 dimensional case can be identified with a real
tree. 

It has been shown by Andersson \cite{An05} that for simplicial
spacetimes generated by deforming a Lorentz cone spacetime with respect to a
finite collection of non-intersecting totally geodesic hypersurfaces, the
above picture holds for the CMC foliation, and the conjecture stated by
Benedetti and Guadagnini \cite{BG01} is valid in this
case.

\subsection*{The Gauss map} 
The Gauss map $\phi: \tM \to \HH^n$ is equivariant with respect to the action
of the isometry group on $\tM \subset \dev(V)$. Therefore, in the case where
the spatial topology is hyperbolic, $\phi$ defines a
Gauss map $\phi: M \to \HH^n /\Gamma$ where $\Gamma$ is the linear part of
the holonomy group of $V$. $\phi: (M,g_\tau) \to (M,\gamma)$ is a
harmonic diffeomorphism, which is isotopic to 
the identity map. 

\subsection*{The reduced Hamiltonian} 
Suppose $V$ is a spatially compact MGHF spacetime, not of translational
type. The mean curvature of the unique global CMC foliation is a time
function and the flow $\tau \mapsto g_\tau$ is the solution to the Einstein
evolution equation with the CMC gauge. By choosing a suitable time gauge $t =
f(\tau)$, one finds that the reduced Hamiltonian is 
$$
\Ham = |\tau|^n \Area(M,g_\tau)
$$
This scale invariant quantity is related to the Yamabe invariant (also known
as the $\sigma$-invariant, in fact \cite[\S 1.2]{An02a}
$$
\inf_{g,K} \frac{n-1}{n} \Ham^{2/n} \geq -\sigma(M)
$$
and in case the Cauchy surface is of hyperbolic type, one conjectures 
$$
\inf_{g,K} \Ham = n^n \Area(M,\gamma) 
$$
where $\gamma$ is a hyperbolic metric on $M$ with sectional curvature
$-1$. Note that $\Area(M,\gamma)$ is a topological invariant. 

By \cite{An02a}, we have 
$$
\lim_{\tau \nearrow 0} \Ham(g,K) = n^n \Area(M,\gamma)
$$
in case the Cauchy surface is of hyperbolic type. 

More generally, when $M$ has a torus factor, the same argument as in
\cite{An02a} gives 
$$
\lim_{\tau \nearrow 0} \Ham(g,K) = 0
$$
since in this case, the torus factor collapses when viewed in the rescaled
geometry. 

\bibliographystyle{amsplain}
\bibliography{data}

\end{document}